\documentclass[11pt]{article}
\usepackage{amsmath,amssymb}
\newcommand{\hess}{{\rm Hess}}
\newcommand{\grad}{{\rm grad}}
\newcommand{\e}{{\rm e}}

\newcommand{\real}{\mathbb{R}}

\newcommand {\cqd}{\begin{flushright}\vskip-25pt$\Box$\end{flushright}}
\newtheorem{myth}{Theorem}[section]
\newtheorem{mylem}{Lemma}[section]
\newtheorem{myprop}{Proposition}[section]

\newtheorem{myrem}{Remark}[section]
\newtheorem{mycoro}{Corollary}[section]
\begin{document}
\title{ Curvature integral estimates for complete hypersurfaces }
\author { Hil\'{a}rio Alencar\thanks{The authors were partially supported by CNPq and FAPERJ, Brazil.},  Walcy
Santos$^*$ {and} {Detang Zhou$^*$ }}
\date {{\em Dedicate to Professor Manfredo do Carmo on the occasion of his 80${th}$ birthday.}}
\maketitle
\vspace*{7pt}
\begin {quote}
{\scriptsize{\bf Abstract.} We consider the integrals of $r$-mean curvatures $S_r$ of a complete hypersurface $M$ in space forms $\mathcal{Q}_c^{n+1}$ which generalize volume $(r=0)$, total mean curvature $(r=1)$, total scalar curvature $(r=2)$ and total curvature $(r=n)$. Among other results we prove that a complete properly immersed hypersurface of a space form with
 $S_r\geq 0$, $S_r\not\equiv 0$ and $S_{r+1}\equiv 0$ for some $r\le n-1$ has
$\int_MS_rdM=\infty.$
}

{\scriptsize {\it Key words}:  mean curvature, space form,
hypersurface, volume estimate.}
\end{quote}
\vspace*{7pt}
\section {Introduction}\indent
Let $M^n$ be a complete orientable  immersed hypersurface of a space
form ${\mathcal Q}_c^{n+1}$ of constant sectional curvature $c$.
Let $A$
be the second fundamental operator of the immersion and let ${\lambda _1,...,\lambda _n }$
be the eigenvalues of $A$.
We define the {\it r-mean curvature} of the immersion at a point $p$ by
$$H_r(p)=\frac{1}{\binom nr}\sum_{i_1<...<i_r}\lambda _{i_1}...\lambda _{i_r}=\frac{1}{\binom nr}S_r(p),$$
where $S_r$ is the $r$ symmetric function of the ${\lambda _1,...,\lambda _n }$, for $1\leq r\leq n$, and define $H_0=1$ and $H_r=0$, for all $r\geq n+1$. In particular, for $r=1$, $H_1=H$ is the mean curvature of the immersion.\\

We define the {\it r-area } of a domain $D\subset M$ by
$${\mathcal A}_r(D)=\int_D S_r(P)\; dM.$$
Then, when $r=0$, ${\mathcal A}_0$ is the volume of $D$.

In this paper we are interested in $r$-areas estimates.
When $r=0$, it is well known that a complete properly immersed minimal hypersurface in $\mathbb{R}^{n+1}$ has at least polynomial volume growth. In fact volume infinity results holds for more general ambient spaces. Precisely we have the following result of K. Frensel \cite {Fr}. \\[10pt]
{\bf Theorem } (\cite{Fr} Theorem 1) {\it Let $M^m$ be a complete,
noncompact manifold and let $x:M^m\rightarrow N^n$ be an isometric
immersion with mean curvature vector field bounded in norm. If $N^n$
has sectional curvature bounded from above and injective radius
bounded from below by a positive constant, then the volume of $M^m$
is infinite.}
\vspace*{10pt}

It is also true that each end of $M$ has infinite volume under the
same conditions(see \cite{CCZ}). These estimates have been used in studying the topology and geometric properties of minimal hypersurfaces and hypersurfaces with constant mean curvature (see for example \cite{Fr}, \cite{CCZ}, \cite{Si}). It would be natural to ask the
following.\\[5pt]

\noindent
{\bf Question }
{\it Let $M^n$ be a complete noncompact manifold and let
$x:M^n\rightarrow N^{n+1}$ be an isometric immersion such that there
is a positive constant $C$ satisfying
\[ |S_{r+1}|\le C S_r,
\]
for some $r=0,1,\cdots, n-1$. Is the $r-$area  of $M^m$ infinite?}

When $r=n$ $S_{r+1}=0$ one can find a negative answer to this question by taking
an example that $M$ is a complete noncompact surface in $\real^3$
with positive Gaussian curvature and the total curvature is finite
by the theorem of Cohn-Vossen. When $r<n$ we obtain $r$-area estimate and give positive answers to the question in some interesting cases.

To state our results we introduce the $r$'th Newton transformation, $P_r:T_pM\rightarrow T_pM$,
which are defined inductively by
$$\begin{array}{ll}
P_0= & I,\\
P_r= &  S_rI-A\circ P_{r-1},\;r>1.
\end{array}$$\\
\noindent
{\bf Theorem A}. (Theorem \ref{thm2}) {\em
Let ${\mathcal Q}_c^{n+1}$ be a Riemannian manifold with constant
sectional curvature $c$ and  $M^n$ a  complete
noncompact properly  immersed hypersurface of ${\mathcal Q}^{n+1}_c$. Assume that
there exists a nonnegative constant $\alpha$ such that
\begin{equation*}(r+1)|S_{r+1}|\le (n-r)\alpha S_r,
\end{equation*}
for some $r\le n-1$. If $P_r$ is positive semi definite, then for
any $q\in M$ such that $S_r(q)\neq 0$ and any $\mu_0>0$ there exists a positive constant $C$
depending on  $\mu_0$, $q$ and $M$ such that
\begin{equation*}
    \int_{\overline{B}_\mu(q)\cap M}S_rdM\ge \int_{\mu_0}^{\mu} C  \e ^{-\alpha
    \tau}d\tau.
\end{equation*}
For the case $c>0$, one needs that $\mu \leq \frac{\pi}{2\kappa}.$}\\[10pt]

As a consequence of this result we obtain the following, which is one of the main results of this article.\\
\noindent {\bf Theorem B}. (Corollary \ref{coro2.2}) {\em Let
${\mathcal Q}_c^{n+1}$ be a Riemannian manifold with constant
sectional curvature $c\leq 0$ and  $M^n$ a  complete noncompact
properly immersed hypersurface of ${\mathcal Q}^{n+1}_c$. Assume
that $S_r\geq 0$, $S_r\not\equiv 0$ and $S_{r+1}\equiv 0$ for some
$r\le n-1$. Then $\int_MS_rdM=\infty.$}
\begin{myrem}
 The cases
when $r$ is even and $r$ is odd are different. If $r$ is odd and
$S_r\leq 0$, we can change the orientation so that $S_r\geq 0$. But
when $r$  is even, $S_r$ is independent of the choice of
orientation. It has been proved by Gromov and Lawson that the
existence of a complete metric with nonpositive scalar curvature
$(r=2)$ implies some topological obstructions, which is called
enlargeable(see Corollary A in \cite{GL}). Enlargeable manifolds
cannot carry metrics of positive scalar curvature.
\end{myrem}
\begin{myrem}This result has been used to study the stable hypersurfaces with constant scalar curvature in Euclidean
spaces in \cite{ASZ}.
\end{myrem}

It is known that the volume estimate of submanifold is related to
the validity of Sobolev inequality. Topping \cite{To} used Sobolev
inequality to get a diameter estimate in terms of the integral of mean curvature.
In the section \ref{sec3}, we get a global estimate of the integral of mean
curvature which is sharp for cylinders. Precisely we prove that\\[10pt]
{\bf Theorem C}. (Theorem \ref{thmmean}) {\em Let $M^m$ be an $m$-dimensional complete noncompact Riemannian manifold
isometrically immersed in $\real^n$. Then there exists a positive
constant $\delta$ depending on  $m$ such that if
\begin{equation*}
    \lim\sup_{r\to +\infty}\frac{V(x,r)}{r^m}<\delta,
\end{equation*}
where $V(x,r)$ denotes the volume of the geodesic ball $B_r(x))$,
then
\begin{equation*}
    \lim\sup_{R\to +\infty}\frac{\int_{B_R(x)}|H|^{m-1}dM}{R}>0.
\end{equation*}
In particular, $\int_M|H|^{m-1}dM=+\infty.$}\\[5pt]

For a complete noncompact surface $M$ with finite total curvature, Cohn-Vossen theorem says that (see Theorem 6 in \cite{CV})

\begin{equation*}
  \int_MKdM\leq 2\pi\chi (M)
\end{equation*}
A special case of Corollary \ref{coro4.1} says that $\int_M|H|dM=+\infty$ if equality holds.\\

The rest of the paper is organized as follows. In Section \ref{sec1} we compute some formulas for distance function and $r$-mean curvature and apply then to main results. The estimate obtained in Section \ref{sec1} can be improved when $r=0$ and this is demonstrated in \ref{sec2}. In Section \ref{sec3} we give the proof of Theorem C.\\[10pt]
\noindent {\bf Acknowledgement.} The authors would like to thank
Professor M. P. do Carmo for  many invaluable comments, suggestions
and encouragements. We would also like to thank M. Dajczer and L.
Florit for  interests and comments.
\section {$r$-area estimate}\label{sec1}
Let $x:M^n\rightarrow N^{n+1}$ be an isometric immersion of a Riemannian manifold into a Riemannian manifold $N$.

In \cite {Re}, Reilly showed that the $P_r$ satisfy the following
\begin{myprop}
(\cite {Re}, p.224, see also \cite{BC}- lemma 2.1) Let
$x:M^n\rightarrow N^{n+1}$ be an isometric immersion between two
Riemannian manifolds and let $A$ be the second fundamental form of
$x$.The $r$'th Newton transformation $P_r$ associated to $A$
satisfies:
\begin{eqnarray}
    \;\mbox{{\rm
trace}}(P_r)&=&(n-r)S_r,\label{eq0.1}  \\
  \;\mbox{{\rm
trace}}(A\circ P_r)&=&(r+1)S_{r+1}. \label{eq0.2}
\end{eqnarray}
 \end{myprop}

 For hypersurfaces with bounded mean curvature, the Laplacian of the
square of the intrinsic distance to a fixed point of $M$ played an
important role in the proof of Frensel's estimate of the volume of
$M$. In the case of  $H_r$ bounded, with $r>1$, we used another
second order differential operator defined on $M$, which seems to be
natural  for this problem. Associated to each Newton transformation $P_r$, if $f:M\rightarrow \mathbb{R}$ is a smooth function, we
define $$L_r(f)=\mbox{{\rm trace}}(P_r\circ\textrm{Hess} \;f).$$

These operators are, in a certain sense, a generalization of the Laplace operator since $L_0(f)=$  $\mbox{{\rm trace}}(\textrm{Hess} \;f)=\Delta f$. They were introduced by Voss \cite{Vo} in connection with variational arguments. In general, these operators are not elliptic and some conditions are necessary to ensure the presence of ellipticity.

We include here some useful facts.
\begin{myprop}\label{2elliptic}(\cite{El}- Lemma 3.10)Let $N^{n+1}$ be an $(n+1)-$dimensional oriented Riemannian manifold and
let $M^{n}$ be a connected $n-$dimensional orientable
Riemannian manifold. Suppose $x:M\rightarrow N$ is
an isometric immersion. If $H_2>0$, then the operator $L_1$ is
elliptic.
\end{myprop}
\begin{myprop}\label{relliptic}(\cite{CR}- Proposition 3.2) Let $N^{n+1}$ be an $(n+1)-$dimensional oriented Riemannian manifold and
let $M^{n}$ be a connected $n-$dimensional orientable Riemannian
manifold (with or without boundary). Suppose $x:M\rightarrow N$ is
an isometric immersion with $H_r>0$ for some $1\leq r\leq n$. If
there exists an interior point $p$ of $M$ such that all the
principle curvatures at $p$ are nonnegative, then for all $1\leq
j\leq r-1,$ the operator $L_j$ is elliptic, and the $j$-mean
curvature $H_j$ is positive.
\end{myprop}

We need the following proposition which is essentially the content of Lemma 1.1 and equation (1.3) of \cite{HL}. We include here with a direct proof.
\begin{myprop}\label{semidef}
Let $M^n\rightarrow N^{n+1}$ be an isometric immersion. Suppose that $S_{r+1}(p)=0$, for some $r$, $0\leq r<n$. Then $P_r$ is semi definite at $p$.
\end{myprop}
{\bf Proof.} Consider $S_r=S_r(\lambda_1,...,\lambda_n)$. Then $\dfrac {\partial S_r}{\partial\lambda_i}$ are the eigenvalues of $P_r$. Let $(\lambda_1^0,...,\lambda_n^0)$ be the principal curvatures of $M$ at $p$. Hence $$S_{r+1}(\lambda_1^0,...,\lambda_n^0)=0.$$
We choose $\epsilon ={\displaystyle \min_{\lambda_i^0\neq 0}\{1, |\lambda_i^0|\}}$. Then, for all $(\varepsilon_1,..., \varepsilon_n)$ with $\varepsilon_i\in(0,\epsilon)$,
$S_{r+1}(\lambda_1^0+\varepsilon_1,...,\lambda_n^0+\varepsilon_n)$ does not change sign. This implies that $\dfrac {\partial S_r}{\partial\lambda_i}\geq 0$ for all $i=1,..,n$ or
$\dfrac {\partial S_r}{\partial\lambda_i}\leq 0$ for all $i=1,..,n$. Thus $P_r$ is semi definite at $p$.\\
\cqd
Let $M^n$ and $N^{n+1}$ be Riemannian manifolds and
$x:M^n\rightarrow N^{n+1}$ an isometric immersion. Henceforth, we
shall tacitly make the usual identification of $X\in T_pM$ with
$dx_p(X)$. In particular, if $F: N\rightarrow \real$ is smooth and
we consider the composition $f=F\circ x$, then we have at $p\in M$,
for every $X\in T_pM$:
$$
\langle \grad_Mf,X\rangle=df(X)=dF(X)=\langle \grad_NF,X\rangle,
$$
where $\grad_M$ and $\grad_N$ denotes the gradient on $M$ and the
gradient on $N$, respectively. So that
\begin{equation}\label{eq0.3}
    \grad_NF=\grad_Mf+(\grad_N F)^\bot,
\end{equation}
where $(\grad F)^\bot$ is perpendicular to $T_pM$.
 Let $F:N\to \mathbb{R}$ be a $C^2$ function and denote $f:M\to
\mathbb{R}$ the function induced by $F$ by restriction, that is
$f=F\circ x$.  We have the following.
\begin{mylem} Let $x:M^{n}\rightarrow N^{n+1}$ an isometric immersion. Let $F:N\rightarrow \real$ a smooth function and considerer $f=F\circ x:M\rightarrow \real$. For  an orthonormal frame  $ \{ e_i\}$ on $M$, we have
\begin{equation}\label{eq1}
    L_rf=\sum_{i=1}^n\mbox{\rm Hess}(F)(e_i, P_r(e_i))+(r+1)S_{r+1}\langle \mbox{\rm grad}_N F, \eta\rangle,
\end{equation}
where $\eta$ denotes the normal vector field of the immersion and
$\mbox{grad}_N$ is the gradient of $N$.
\label{lemma2.1}\end{mylem} {\bf Proof}. Let $\nabla$ and
$\overline{\nabla}$ the connection of $M$ and $N$, respectively. If
$\alpha$ denotes the second fundamental form of the immersion,
Gauss' equation and equations (\ref{eq0.2}) and (\ref{eq0.3}) imply
that

\begin{equation*}
    \begin{split}
       L_rf &= \sum_{i=1}^n\langle\nabla_{e_i}(\mbox{grad}_M
         f),P_r(e_i)\rangle\\
         &=\sum_{i=1}^n\langle\overline{\nabla}_{e_i}(\mbox{grad}_Mf)-[\overline{\nabla}_{e_i}(\mbox{grad}_M f)-\nabla_{e_i}(\mbox{grad}_M
         f)],P_r(e_i)\rangle\\
         &=\sum_{i=1}^n\langle\overline{\nabla}_{e_i}(\mbox{grad}_Mf)-\alpha (e_i, \mbox{grad}_M f),P_r(e_i)\rangle\\
         \end{split}
         \end{equation*}
         \begin{equation*}
         \begin{split}
         &=\sum_{i=1}^n\langle\overline{\nabla}_{e_i}(\mbox{grad}_Mf),P_r(e_i)\rangle\\
         &=\sum_{i=1}^n\langle\overline{\nabla}_{e_i}(\grad_NF-(\grad_NF)^\bot),P_r(e_i)\rangle\\
        &=\sum_{i=1}^n\langle\overline{\nabla}_{e_i}\grad_NF,P_r(e_i)\rangle
        -\sum_{i=1}^n\langle\overline{\nabla}_{e_i}(\grad_NF)^\bot,P_r(e_i)\rangle\\
        &=\sum_{i=1}^n\hess(F)(e_i,P_r(e_i))
        -\sum_{i=1}^n\langle\overline{\nabla}_{e_i}(\langle \grad_NF, \eta\rangle \eta),P_r(e_i)\rangle\\
         \end{split}\end{equation*}
        \begin{equation*}\begin{split}
        &=\sum_{i=1}^n\hess(F)(e_i,P_r(e_i))
        -\sum_{i=1}^n\left\langle\langle \grad_NF, \eta\rangle \overline{\nabla}_{e_i}\eta,P_r(e_i)\right\rangle\\
        &=\sum_{i=1}^n\hess(F)(e_i,P_r(e_i))
        -\langle \grad_NF, \eta\rangle\sum_{i=1}^n\left\langle -A(e_i),P_r(e_i)\right\rangle\\
       \end{split}
\end{equation*}
\begin{equation*}
    \begin{split}
    &=\sum_{i=1}^n\hess(F)(e_i,P_r(e_i))
        +\langle \grad_NF, \eta\rangle\sum_{i=1}^n\left\langle e_i,AP_r(e_i)\right\rangle\\
        &=\sum_{i=1}^n\hess(F)(e_i,P_r(e_i))
        +\langle \grad_NF, \eta\rangle {\rm trace (AP_r)}\\
        &=\sum_{i=1}^n\textrm{Hess}(F)(e_i,
P_r(e_i))+(r+1)S_{r+1}\langle\grad_N F, \eta\rangle.
 \end{split}
\end{equation*}\\\cqd

Let $c\in \real$. Define the function
$c_\kappa(t)=\displaystyle{\int_0^t}s_\kappa(t)dt$ where
\begin{equation}\label{eq2}
    s_\kappa(t)=\left\{
                  \begin{array}{ll}
                    \frac{\sin \kappa t}{\kappa}, & \hbox{ if $c=\kappa^2$;} \\
                    t, & \hbox{if $c=0$;} \\
                    \frac{\sinh \kappa t}{\kappa}, & \hbox{if $c=-\kappa^2$.}
                  \end{array}
                \right.
\end{equation}
 If $\rho $ denotes the distance function to the point $Q$ in
$N^{n+1}$, let $F:N^{n+1}\rightarrow\mathbb{R}$ given by
$F(p)=c_\kappa(\rho(p))$. Therefore the lemma \ref{lemma2.1} with
$f=F\circ x,$ where $F=c_\kappa \circ \rho$ implies
\begin{mycoro}\label{cor1}
Let $M$ be an immersed hypersurface in $N^{n+1}$  and let $\kappa
\in \real$. Then, for all $p\in M$,
\begin{equation}
L_r\left(c_\kappa(\rho(p))\right)=(n-r)s'_\kappa(\rho(p))S_r+(r+1)S_{r+1}s_{\kappa}(\rho(p))\langle\mbox{\rm
grad}_{N} \rho(p),\eta\rangle.\label{eqcoro2.1}
\end{equation}
In particular, when $c=0,$
\[
\frac
12L_r\left(\rho^2(p)\right)=(n-r)S_r+(r+1)S_{r+1}\rho(p)\langle\mbox{\rm
grad}_{N} \rho(p),\eta\rangle.
\]
\end{mycoro}
{\bf Proof.} First observe that
\begin{equation}\label{eq0.4}
    \hess F(X,Y)=s_\kappa(\rho)\langle X,Y \rangle,
\end{equation}
where $X,Y \in T_{x(p)}\mathcal{Q}$. In fact,
\begin{equation*}\begin{split}
\hess F(X,Y)&=\hess (c_\kappa (\rho))\\
&=\left \langle \overline{\nabla}_X\grad_N (c_\kappa (\rho )), Y
\right \rangle\\
&=\left \langle \overline{\nabla}_X s_\kappa (\rho)\grad_N  \rho , Y
\right \rangle\\
&=s_\kappa (\rho) \hess \rho (X,Y) + s^\prime_\kappa \left \langle
\langle \grad_N\rho, X\rangle \grad_N\rho , Y
\right \rangle.\\
\end{split}
\end{equation*}
On the other hand, see \cite{AF}, p.6,
$$\hess \rho(X,Y)=\left\langle
\overline{\nabla}_X\grad_N\rho,Y\right\rangle=\frac{s^\prime_\kappa(\rho)}{s_\kappa
(\rho)}\left [\langle X,Y \rangle -\langle \grad_N \rho, X\rangle
\langle \grad_N \rho, Y\rangle\right ].$$ This concludes the proof
of (\ref{eq0.4}). Now, by using equation (\ref{eq0.3}), we have
\begin{equation*}\begin{split}
L_rf &= \sum_{i=1}^ns_\kappa^\prime(\rho) \langle e_i,P_r(e_i)\rangle -(r+1)S_{r+1}\langle \grad_N(c_\kappa\circ \rho), \eta\rangle\\
&=s^\prime_\kappa(\rho){\rm
trace}P_r-(r+1)S_{r+1}s_\kappa(\rho)\langle\grad_N\rho,\eta\rangle.
\end{split}\end{equation*}
Finally, by using equation (\ref{eq0.1}), we conclude the proof of
equation (\ref{eqcoro2.1}). The case $c=0$ follows immediately.
\\ \cqd

Let ${\mathcal Q}_c^{n+1}$ be a Riemannian manifold with constant
sectional curvature $c$ and let $x:M\rightarrow {\mathcal Q}_c^{n+1}$ an isometric immersion. It follows from Codazzi equation (see
Rosenberg \cite {Ro}, p.225) that $L_r$ is a divergent form
operator, that is,
$$L_r(f)=\textrm{div}_M(P_r\nabla f),$$
for all smooth function $f:M\rightarrow \mathbb{R}$. Denote by $B_r(Q)$ the geodesic ball of
${\mathcal Q}_c^{n+1}$ with radius $r$, and center $Q\in {\mathcal Q}_c^{n+1}$  and by
$\overline{B}_r(Q)$ its closure. We will use the following
proposition to prove our results.
\begin{myprop}\label{prop3} Let ${\mathcal Q}_c^{n+1}$ be a Riemannian manifold with constant
sectional curvature $c$ and  let $x:M^n\rightarrow {\mathcal
Q}_c^{n+1}$ an  isometric immersion. For $Q\in {\mathcal
Q}_c^{n+1}$, we denote by $\rho(x)$ the distance to the point $Q\in
{\mathcal Q}_c^{n+1}$ and $\rho(x(p))$, $p\in M$ its restriction to
$M$. If for some $r\leq n-1$, $S_r\ge 0$, then
\begin {equation}\begin{split}
\int_{\partial D}s_\kappa(\rho(q))\langle P_r(\mbox{\rm
grad}_\mathcal{Q} \rho(q)^\top), \nu\rangle dA&\geq
(n-r)\int_{D}(s_\kappa'(\rho(q))S_r-\frac{r+1}{n-r}|S_{r+1}|
s_\kappa(\rho(q))) dM, \label {eq3.5}
\end{split}
\end {equation}
where $q=x(p)$, $\nu$ is the conormal vector of $D$ and $D\subset M$ is a bounded domain with nonempty boundary
$\partial D$. In the case $c>0$, we also request that $D\subset \overline{B}_{\frac{\pi}{2\kappa}}(Q)$.
\end{myprop}
{\bf Proof}. By using (\ref{eqcoro2.1}), and since
$|\grad_\mathcal{Q}\rho(x(p))|\leq 1$,
$s^\prime_\kappa(\rho(x(p))\geq 0$, we have

\[  L_r (c_\kappa(\rho(x)))\geq (n-r)[s_\kappa'(\rho)
S_r-\frac{r+1}{n-r}|S_{r+1}|s_\kappa(\rho)]. \]
  Integrating this inequality,
we get
\begin {equation}
\int_{D} L_r (c_\kappa(\rho(x)))dM\geq
(n-r)\int_{D}[s_\kappa'(\rho(x))S_r-\frac{r+1}{n-r}|S_{r+1}|
s_\kappa(\rho(x))] dM.
\end {equation}
On the other hand, using the divergence theorem, we have that
\begin{equation*}
\begin{split}
\int_{D}L_r (c_\kappa(\rho(x)))dM &=\int_D\textrm{div}
P_r(\grad_M(c_\kappa(\rho(x(p)))))\\
&=\int_D\textrm{div} (s_\kappa\rho(x(p))P_r(\grad_\mathcal{Q}\rho)^\top)\\
 &=\int_{\partial D} s_\kappa(\rho(x))\langle P_r((\mbox{\rm
grad}_\mathcal{Q}\rho)^\top), \nu\rangle dA,
\end{split}\end{equation*}
 where $\nu$ denotes the outward unit normal vector field on $\partial D$. Therefore, if $q=x(p)$
\begin{equation*}
    \begin{split}
        \int_{\partial D}
s_\kappa(\rho(q))\langle P_r((\mbox{\rm
grad}_\mathcal{Q}\rho(q))^\top), \nu\rangle dA\ge &
(n-r)\int_{D}[s_\kappa'(\rho(x))S_r-\frac{r+1}{n-r}|S_{r+1}|
s_\kappa(\rho(x))] dM,
     \end{split}
\end{equation*}
and the proposition is proved.\\\cqd

Observe that the above Proposition is valid for a more general class of domains. For instance it is valid in the setting of Gauss-Green Theorem (see \cite{Fe}, p.478). In particular, if we take $D$ to be the intersection of the extrinsic ball with $M$
i.e. $D=\overline{ B}_{\mu}\cap M$ in Proposition \ref{prop3}, we have the
 following
 \begin{myth}\label{thm2}
Let ${\mathcal Q}_c^{n+1}$ be a Riemannian manifold with constant
sectional curvature $c$ and  $M^n$ a  complete
noncompact properly  immersed hypersurface of ${\mathcal Q}^{n+1}_c$. Assume that
there exists a nonnegative constant $\alpha$ such that
\begin{equation}(r+1)|S_{r+1}|\le (n-r)\alpha S_r,\label{eq4.05}
\end{equation}
for some $r\le n-1$. If $P_r$ is semi-positive definite, then for
any $q\in M$ such that $S_r(q)\neq 0$ and any $\mu_0>0$ there exists
a positive constant $C$ depending on  $\mu_0$, $q$ and $M$ such that
\begin{equation*}
    \int_{\overline{B}_\mu(q)\cap M}S_rdM\ge \int_{\mu_0}^{\mu} C  \e ^{-\alpha
    \tau}d\tau.
\end{equation*}
For the case $c>0$, one needs that $\mu \leq \frac{\pi}{2\kappa}.$
\end{myth}
\noindent
{\bf Proof.} We use the notation introduced in Proposition \ref{prop3}.
 Take  $D=D_\tau= \overline{B}_\tau(q)\cap M$. Since the immersion is proper, we have that $\partial D_\tau \neq \emptyset$, for all $\tau>0$. Thus, by using (\ref{eq4.05}) in equation (\ref{eq3.5}), we obtain that
\begin {equation}\begin{split}
\int_{\partial D_\tau}{s_\kappa(\rho (x))}&\langle P_r(\mbox{\rm grad}_M \rho ),\nu\rangle dA\geq
(n-r)\int_{D_\tau}(s_\kappa'(\rho (x))-\alpha
s_\kappa(\rho (x))) S_rdM\\
&= (n-r)\int_0^\mu\int_{\partial
D_\tau}\frac{(s_\kappa'(\rho (x))-\alpha
s_\kappa(\rho (x))}{s_\kappa(\rho (x))}s_\kappa(\rho (x))|\mbox{\rm grad}_M \rho |^{-1}S_r
dAd\tau, \label {eq4.1}
\end{split}
\end {equation}
where  we have used the co-area formula (see  \cite {Be}, p. 80). Observe that the conormal vector $\nu$ to $\partial D$ is
parallel to $\mbox{\rm grad}_M \rho $. This fact and that $P_r$ is semi-positive
definite, implies that $$\langle P_r(\mbox{\rm grad}_M \rho ),\nu\rangle\le tr(P_r)
|\mbox{\rm grad}_M \rho |= (n-r)S_r| \mbox{\rm grad}_M \rho |.$$
Using the above equation and the fact that along $\partial D_\tau$, $\rho (x)=\tau$,we get
\begin {equation}\begin{split}
\int_{\partial D_\mu} s_\kappa(\rho (x)) |\mbox{\rm grad}_M \rho |S_rdA
&\geq\int_0^\mu\frac{s_\kappa'(\tau)-\alpha
s_\kappa(\tau)}{s_\kappa(\tau)}\int_{\partial
D_\tau}s_\kappa(\rho (x))|\mbox{\rm grad}_M \rho |^{-1} S_rdAd\tau.
\end{split}\label{eq4.2}
\end {equation}
Now we define
\[\varphi(\tau)=\int_{\partial
D_\tau}s_\kappa(\rho (x))|\mbox{\rm grad}_M|^{-1}S_r dA.
\]
The equation (\ref{eq4.2}) implies
\[\varphi(\mu)\ge \int_0^\mu\frac{s_\kappa'(\tau)-\alpha
s_\kappa(\tau)}{s_\kappa(\tau)}\varphi(\tau)d\tau.
\]
If we write
\[\phi(\mu)=\int_0^\mu\frac{s_\kappa'(\tau)-\alpha
s_\kappa(\tau)}{s_\kappa(\tau)}\varphi(\tau)d\tau,
\]
we have
\[\begin{split}
 \phi'(\mu) & \ge \frac{s_\kappa'(\mu)-\alpha
s_\kappa(\mu)}{s_\kappa}\phi(\mu).
\end{split}
 \]
Thus, by integrating from $\mu_0 >0$ to $\mu$, the above differential inequality arises
$$\ln \frac {\phi(\mu)}{\phi(\mu_0)}\geq \ln (\frac {s_\kappa(\mu)}{s_\kappa(\mu_0)}) - \alpha (\mu-\varepsilon)=
\ln [  (\frac {s_\kappa(\mu)}{s_\kappa(\mu_0)})\e
^{-\alpha(\mu-\mu_0)}].$$ Hence,

\[\phi(\mu)\ge \frac{\phi(\mu_0)}{s_\kappa(\mu_0)}s_\kappa(\mu)  \e ^{-\alpha \mu}.
\]

 Define
$$f(\mu)=\int_{D_\mu(q)} S_rdM.$$
Again by the co-area formula, it follows that
$$f(\mu)=\int_0^\mu(\int_{\partial D_\tau(q)}|\mbox{\rm grad}_M \rho |^{-1}S_r dA)d\tau.$$
Since
$$f^\prime (\mu)=\int_{\partial D_\mu(q)}|\mbox{\rm grad}_M \rho |^{-1}S_r dA=\frac 1{s_\kappa(\mu)}\varphi(\mu)\ge  \frac{\phi(\mu_0)}{s_\kappa(\mu_0)}  \e ^{-\alpha \mu},$$
then for  $\mu>\mu_0$,
\begin{equation*}
    f(\mu)\ge \int_{\mu_0}^{\mu} \frac{\phi(\mu_0)}{s_\kappa(\mu_0)}  \e ^{-\alpha
    \tau}d\tau.
\end{equation*}\\
\cqd
 \begin{mycoro}\label{coro2.2}
Let ${\mathcal Q}_c^{n+1}$ be a Riemannian manifold with constant
sectional curvature $c\leq 0$ and  $M^n$ a  complete
noncompact properly  immersed hypersurface of ${\mathcal Q}^{n+1}_c$. Assume that $S_r\geq 0$, $S_r\not\equiv 0$ and $S_{r+1}\equiv 0$ for some $r\le n-1$. Then
$\int_MS_rdM=\infty.$
\end{mycoro}
{\bf Proof}. Since the immersion is proper, we have $\partial (M\cap \overline{ B}_\mu(q))$ is nonempty for all $\mu>\mu_0$. By using Proposition \ref{semidef}, since $S_{r+1}=0$,  we have that $P_r$ is semi-definite. Now, the condition $S_r\geq 0$ implies that $P_r$ is positive semi-definite. Therefore, using Theorem \ref{thm2}, with $\alpha=0$,
for all $\mu>\mu_0$,
$$  \int_{\overline{B}_\mu\cap M}S_rdM\ge \int_{\mu_0}^{\mu} C  \e ^{-\alpha
    \tau}d\tau=C(\mu-\mu_0).$$
Then
$$\int_MS_rdM=\infty.$$\\ \cqd
\begin{myrem}
When $r$ is odd, the condition $S_r\geq 0$ can be obtained by choosing the right orientation.
\end{myrem}
The condition of semi-positiveness of $P_2$  is satisfied when $M$
is hypersurface immersed in $\real^{n+1}$ with $S_3=0$(which is
called 2-minimal) and $S_2>0$. In fact, in this case $P_2$ is
positive definite, since $L_2$ is elliptic (see Proposition
\ref{2elliptic}). So we have

\begin{mycoro}\label{coro1.3}
Let $M^n$ be a complete 2-minimal noncompact properly immersed
hypersurface of $\real^{n+1}$ with   nonnegative scalar curvature.
Then either the  scalar curvature is  zero or the total scalar
curvature is infinite.
\end{mycoro}
\begin{myrem}When $n=3$ the corollary can be proved using
Theorem III in \cite{HN} without the assumption that the immersion
is proper. In this case, $M^n$ has to be cylinder and the conclusion
of the above Corollary follows immediately.
\end{myrem}
\begin{myrem}
 The condition of semi-positiveness of $P_r$ is also satisfied when $M$ is a hypersurface  in
$\real^{n+1}$ with nonnegative sectional or positive Ricci
curvature. Indeed when $Ric_M>0$, for each point in $M$ the
principal curvatures can be arranged as
$\lambda_1\le\lambda_2\cdots\le\lambda_i<0<\lambda_{i+1}\le\cdots\le\lambda_n$.
The positivity of Ricci curvature implies
\[Ric_M(e_j)=\lambda_j(\sum_{k\ne j}\lambda_k ) >0,\quad \forall j=1,...,n.\]
We can see that if $i\neq 1$ and $i\neq n-1$, and
\begin{equation}
         \sum_{k\ne j}\lambda_k <0, \textrm{ when }j\le i,\label{star}
\end{equation}
\begin{equation}
\sum_{k\ne j}\lambda_k  >0, \textrm{ when }j> i.\label{star1}
\end{equation}
From (\ref{star}) we have for $j\leq i$,

$$
\left( \sum_{k=1}^i \lambda_k-\lambda_i\right)+\sum_{k=i+1}^n\lambda_k<0,
$$
Thus
$$
-\sum_{k=1}^i\lambda_k>-\sum_{k=1}^i\lambda_k+\lambda_j>\sum_{k=i+1}^n\lambda_k.
$$
On the other hand, (\ref{star1}), with $j>i$, implies that
$$\sum_{k=1}^i \lambda_k+\sum_{k=i+1}^n\lambda_k-\lambda_j
>0$$
Hence,
$$
-\sum_{k=1}^i\lambda_k<-\sum_{k=1}^i\lambda_k+\lambda_j<\sum_{k=i+1}^n\lambda_k,
$$
which is a contradiction. One can easily see that the cases $i=1$
and $i=n-1$ also can not occur. Thus, all $\lambda_i$ has the same
sign (we are indebted to  F. Fontenele for this observation).  So we
can choose an orientation such that $P_r$ is positive definite and
$S_r>0$.
\end{myrem}
 Thus we have.
\begin{mycoro}

Let $M^n$ be a complete noncompact properly immersed hypersurface of
$\real^{n+1}$ with positive Ricci curvature. Assume that  there
exists a positive constant $\alpha$ such that
\[(r+1)|S_{r+1}|\le (n-r)\alpha S_r,
\]
for some $r\le n-1$. Then for any $q\in M$ and any $\mu_0>0$ there
exists a positive constant $C$ depending on  $\mu_0$, $Q$ and $M$
such that
\begin{equation*}
    \int_{\bar B(\mu)\cap M}S_rdM\ge \int_{\mu_0}^{\mu} C  \e ^{-\alpha
    \tau}d\tau,
\end{equation*}
where  $\overline{B}_{\mu}(q)$ is the geodesic ball in ${\real}^{n+1}$
centered at $q$.

\end{mycoro}

The following is a direct consequence of Theorem \ref{thm2} and Proposition \ref{relliptic}.

\begin{mycoro}
Let $M^n$ be a complete noncompact properly immersed hypersurface of
$\mathcal{Q}_c^{n+1}$. Assume that $S_r$ is positive and there
exists a positive constant $\alpha$ such that
\[(r+1)|S_{r+1}|\le (n-r)\alpha S_r,
\]
for some $r\le n-1$. If there exists a point such that all principal
curvatures are nonnegative, then for
any $q\in M$ and any $\mu_0>0$ there exists a positive constant $C$
depending on  $\mu_0$, $Q$ and $M$
such
that
\begin{equation*}
    \int_{\bar B(\mu)\cap M}S_rdM\ge \int_{\mu_0}^{\mu} C  \e ^{-\alpha
    \tau}d\tau,
\end{equation*}
where  $\overline{B}_{\mu}(q)$ is the geodesic ball in ${\mathcal Q}_c^{n+1}$
centered at $q$. For the case $c>0$, one needs that $\mu \leq \frac{\pi}{2\kappa}.$
\end{mycoro}

\section {Volume estimates in general manifolds}\label{sec2}

In this section we consider $N^{n+p}$ with sectional curvature bounded from above by a constant $c$.
Let $M^n$ be a submanifold  isometrically immersed in $N= N^{n+p}$.

Let $F:
N \to \mathbb{R}$ be a $C^2$ function and denote $f:M\to \mathbb{R}$
the function induced by $F$ by restriction. Essentially in the same way we prove Lemma \ref{lemma2.1}, we obtain
\[\Delta f=\sum_{i=1}^n\textrm{Hess} F(e_i,e_i)+n\langle \mbox{\rm grad}_N
F,\mathbf{H}\rangle,
\]
where $\{ e_1, e_2, \cdots, e_n\}$ is an orthonormal frame along $M$
and $\mathbf{H}$ is the mean curvature vector.
Similar to Proposition \ref{prop3}, we have
\begin{myprop}\label{prop4} Let $N$ be a Riemannian manifold with
sectional curvature bounded above by $c$ and  $M^n$ an  immersed
connected submanifold of $N$. We denote by $\rho(x)$ the
distance between $x$ and $Q\in N^{n+p}$ and $\rho (x)$ the induced function of
$\rho$ by restriction. Then
\begin {equation}\begin{split}
\int_{\partial D}s_\kappa(\rho (x))\langle \mbox{\rm grad}_M \rho , \nu\rangle dA&\geq
n\int_{D}(s_\kappa'(\rho (x))-|\mathbf{H}| s_\kappa(\rho (x))) dM, \label
{eq3.52}
\end{split}
\end {equation}
where $\kappa =\sqrt{|c|}$, $\nu$ is the conormal vector of $D$ and $D\subset M$ is a bounded domain with nonempty boundary
$\partial D$ and $D\bigcap C_N(Q)=\emptyset$, where $C_N(Q)$ is the cut locus of the point $Q$ in $N$.

\end{myprop}
{\bf Proof}.
Let $V=s_\kappa(\rho)\mbox{\rm grad}_N \rho$ and $V^\top$ the
orthogonal projection of $V$ into the tangent space of $M$. Then we
have $V^\top=s_\kappa(\rho ) \mbox{\rm grad}_M \rho $, where $\rho (x)$ is the induced
function of $\rho$ to $M$ by restriction.
Thus, Lemma 2.5 of \cite {JK}, p. 713, implies that when $\rho<\mathrm{inj}_N(Q)$,
\begin{equation}\label{eq3.51}
    \textrm{Hess} F(X,X)\ge s^\prime_{\kappa}(\rho)\langle X,X\rangle.
\end{equation}

Then
$$\langle \overline{\nabla}_{e_i}V, e_i\rangle\ge s_\kappa'(\rho),$$
for all $\rho$ when  $c\le0$, and $\rho \leq \frac{\pi}{\kappa}, $
when  $c>0$. We have that
\[  \Delta (c_\kappa(\rho(x)))\ge n[s_\kappa'(\rho)- s_\kappa(\rho)|\mathbf{H}|].
\]
  Integrating this inequality and
applying Stokes' formula, we get

\begin{equation*}
        \int_{\partial D}
s_\kappa\langle (\mbox{\rm grad}_N \rho)^\top, \nu\rangle dA\ge
n\int_{D}[s_\kappa'(\rho (x))- s_\kappa(\rho (x))|\mathbf{H}|] dM.
\end{equation*}
The proposition is proved.

\cqd

Similar to Proposition \ref{prop3}, the above result is valid in a more general setting, as extrinsic geodesic balls. Using this fact, we  get
\begin{myth}Let $M$ be a Riemannian manifold  isometrically immersed in a geodesic
ball $\bar B(O, \rho_0)\subset N^{n+p}$ with codimension $p$. Assume
that the sectional curvature of $N$ in $\bar B(O, \rho_0)$ is
bounded above by $c$ and  there exists a positive constant $\alpha$
such that
\[|\mathbf{H}|\le \alpha.
\]
 Then
\[\mathrm{vol}( B_\mu(Q)) \geq
n \omega_n\int_0^\mu s_\kappa(s)^{n-1}  \e ^{-n\alpha s}ds,\]
 where $\kappa =\sqrt{|c|}$, $\omega _n$ is the volume of the
unit ball in $\real ^n$ and $B_\mu(q)$ is the intrinsic geodesic
ball in $M$ with center $q\in M$ and radius $\mu< \mathrm{inj}_N(q)$.
\end{myth}

{\bf Proof.} Taking  $D=B_\tau(q)$ in Proposition \ref{prop4}, then
\[\langle \mbox{\rm grad}_M \rho , \nu\rangle\le |\mbox{\rm grad}_M \rho |,
\]
we have
\begin {equation*}\begin{split}
\int_{\partial B_\tau(q)}\frac{s_\kappa(\rho (x))}n |\mbox{\rm grad}_M\rho |dA&\geq
\int_{B_\tau(q)}(s_\kappa'(\rho (x))-\alpha
s_\kappa(\rho (x))) dM\end{split}\end{equation*}\begin{equation}\begin{split}
&= \int_0^\mu\int_{\partial B_\tau(q)}\frac{(s_\kappa'(\rho (x))-\alpha
s_\kappa(\rho (x))}{s_\kappa(\rho (x))}s_\kappa(\rho (x))|\mbox{\rm grad}_M\rho |^{-1}
dAd\tau, \label {eq4}
\end{split}
\end {equation}
where  we have used the co-area formula (see  \cite {Be}, p. 80).
Since the intrinsic distance is not less than the extrinsic one and
\[
\left(\frac{s_\kappa'}{s_\kappa}\right)'\le 0,
\]
then
\begin {equation}\label{eq6}\begin{split}
\frac{1}n\int_{\partial B_\mu(q)} s_\kappa(\rho (x)) |\mbox{\rm grad}_M\rho |dA
&\geq\int_0^\mu\frac{s_\kappa'(\tau)-\alpha
s_\kappa(\tau)}{s_\kappa(\tau)}\int_{\partial
B_\tau(q)}s_\kappa(\rho (x))|\mbox{\rm grad}_M \rho |^{-1} dAd\tau.
\end{split}
\end {equation}
Now we define
\[\varphi(\tau)=\int_{\partial
B_\tau(q)}s_\kappa(\rho (x))|\mbox{\rm grad}_M \rho |^{-1} dA.
\]
Equation (\ref{eq6}) implies
\[\frac1n\varphi(\mu)\ge \int_0^\mu\frac{s_\kappa'(\tau)-\alpha
s_\kappa(\tau)}{s_\kappa(\tau)}\varphi(\tau)d\tau.
\]
If we write
\[\phi(\mu)=\int_0^\mu\frac{s_\kappa'(\tau)-\alpha
s_\kappa(\tau)}{s_\kappa(\tau)}\varphi(\tau)d\tau,
\]
we have
\[\begin{split}
 \phi'(\mu) & \ge \frac{n(s_\kappa'(\mu)-\alpha
s_\kappa(\mu))}{s_\kappa(\mu)}\phi(\mu).
\end{split}
 \]
Thus, by integrating from $\varepsilon >0$ to $\mu$, with $\mu\leq \min \{ inj_N
(q), \frac { \pi}{2\kappa}\}$, when $c>0$, the above differential inequality
arises
$$\frac 1n\ln \frac {\phi(\mu)}{\phi(\varepsilon)}\geq \ln (\frac {s_\kappa(\mu)}{\varepsilon}) - \alpha (\mu-\varepsilon)=
\ln [  (\frac {s_\kappa(\mu)}{\varepsilon})\e
^{-\alpha(\mu-\varepsilon)}].$$ Hence,
\begin {equation}
 \frac {\phi (\mu) }{\phi(\varepsilon)}\geq   \left[(\frac {s_\kappa(\mu)}{\varepsilon})\e ^{-\alpha(\mu-\varepsilon)}\right]^n.
\label {eq5}
\end {equation}
Observe that by the mean value theorem,
$$\lim_{\varepsilon\rightarrow 0} \frac {\phi(\varepsilon)}{\varepsilon
^n}=\omega_n.
$$  Then
\[\phi(\mu)\ge \omega_ns_\kappa(\mu)^n  \e ^{-n\alpha s}.
\]

 Now, define
$$f(\mu)=\int_{B_\mu(q)} dM=vol(B_\mu(q)).$$
Again by the co-area formula, it follows that
$$f(\mu)=\int_0^\mu(\int_{\partial B_\tau(q)}|\mbox{\rm grad}_M \rho |^{-1} dA)d\tau.$$
Hence
$$f^\prime (\mu)=\int_{\partial B_s(q)}|\mbox{\rm grad}_M \rho |^{-1} dA .$$
This equality and $|\mbox{\rm grad}_M\rho |\le 1$,  with equation (\ref {eq4})
imply that
$$\frac{s_\kappa(\mu)}n{f^\prime (\mu) }\geq \int_{\partial B_\tau(q)}\frac{s_\kappa(\rho (x))}n |\mbox{\rm grad}_M\rho |dA\geq \int_0^\mu(s_\kappa'(\tau)-\alpha
s_\kappa(\tau))f^\prime (\tau)d\tau.$$
 Since \[ f'(\mu)\ge
\frac{\varphi(\mu)}{s_\kappa(\mu)},\] then
\[f(\mu)\ge \int_0^\mu \omega_nns_\kappa(\tau)^{n-1}  \e ^{-n\alpha \tau}d\tau,\]
which concludes the proof.
\cqd

From the theorem we have an immediate corollary.
 \begin {mycoro}\label {coro1}
(i) Let   $M^n$ an  immersed minimal hypersurface of the Euclidean
space ${\mathbb{R}}^{n+p}$. Then
\[vol( B_\mu(q)) \geq
\omega_n\mu^n. \]
where  $\omega _n$ is the volume of the
unit ball in $\real ^n$ and $B_\mu(q)$ is the intrinsic geodesic
ball in $M$ with center $q\in M$.\\
(ii) Let   $M^n$ an  immersed  hypersurface of the hyperbolic space
${\mathbb{H}}^{n+p}(-1)$. Assume there exists a positive constant
$\alpha$ such that
\[|H|\le \alpha<\frac{n-1}n.
\]Then exist a constant $C>0$ such that if $\mu \geq 1$,
\[vol( B_\mu(q)) \geq
Ce^{(n-1-n\alpha)\mu}, \] where  $B_\mu(q)$ is the intrinsic geodesic
ball in $M$ with center $q\in M$.

\end {mycoro}

\section {Mean curvature  integral}\label{sec3}
In this section, inspired by a recent work of Topping \cite{To}, we prove a type of mean curvature integral estimate for complete submanifold in a Euclidean space $\real^n$ and we apply it to surfaces in $\real^n$
\begin{myth}\label{thmmean}Let $M^m$ be an $m$-dimensional complete noncompact Riemannian manifold
isometrically immersed in $\real^n$. Then there exists a positive
constant $\delta$ depending on  $m$ such that if
\begin{equation}\label{eq20}
    \lim\sup_{r\to +\infty}\frac{V(x,r)}{r^m}<\delta,
\end{equation}
where $V(x,r)$ denotes the volume of the geodesic ball $B_r(x))$,
then
\begin{equation}\label{eq21}
    \lim\sup_{R\to +\infty}\frac{\int_{B_R(x)}|H|^{m-1}dM}{R}>0.
\end{equation}
In particular, $\int_M|H|^{m-1}dM=+\infty.$
\end{myth}
We need the following lemma of Topping \cite{To}.
\begin{mylem}\label{lemTopping}(\cite{To}, Lemma 1.2) Let $M^m$ be an $m$-dimensional complete Riemannian manifold
isometrically immersed in $\real^n$. Then exists a positive
constant $\delta$ depending on  $m$ such that for any $x\in M$ and
$R>0$, at least one of the following is true:\\
\hspace*{10pt}$(i)\quad{\displaystyle \sup_{r\in(0,
    R]}r^{-\frac1{m-1}}[V(x,r)]^{-\frac{m-2}{m-1}}\int_{B(x,r)}|H|^{m-1}dM>\delta,} $\\
\hspace*{8pt} $(ii)
\quad{\displaystyle
    \inf_{r\in(0,
    R]}\frac{V(x,r)}{r^m}>\delta.}
$
\end{mylem}
{\bf Proof of Theorem \ref{thmmean}.} We can choose $L$ large enough so that $V(z,L)\le \delta L^m$ for
all $z\in M$. Then from Lemma \ref{lemTopping}, we have
\begin{equation*}
    \sup_{r\in(0,
    L]}r^{-\frac1{m-1}}[V(z,r)]^{-\frac{m-2}{m-1}}\int_{B_r(z)}|H|^{m-1}dM>\delta.
 \end{equation*}
Since
\begin{equation*}
    \int_{B_r(z)}|H|dm\le
    \left(\int_{B_r(z)}|H|^{m-1}dM\right)^\frac1{m-1}\cdot
    \left(V(z,r)\right)^\frac{m-2}{m-1},
\end{equation*}
for any $z\in M$, there exists a $r(z)\in(0, R]$ such that
\begin{equation*}
    \int_{B_r(z)}|H|^{m-1}dM> \delta^{m-1}r(z).
\end{equation*}
Fix a point $o\in M$, we can find a ray $\gamma:[0,+\infty)\to M$
parameterized by arclength. For any fixed $R>0$,
\[\gamma([0,R])\subset\bigcup_{t\in[0,R]}B_{r(\gamma(t))}(\gamma(t)).
\]
From a covering argument used in Theorem 1.1 of \cite{To}, we can find  an at most countable sequence
$t_1, t_2, \cdots, t_q,\cdots \in [0,R]$ such that
$\sum_{i}r(\gamma(t_i))\ge \frac 14R$ and when $i\ne j$
\begin{equation*}
    B_{r(\gamma(t_i))}(\gamma(t_i))\bigcap B_{r(\gamma(t_j))}(\gamma(t_j))=\emptyset.
\end{equation*}
Then
\begin{equation*}
    \begin{split}
      \int_{B_{2R}(o)}|H|^{m-1}dM & \ge \sum_i\int_{B_{r(\gamma(t_i))}(\gamma(t_i))}|H|^{m-1}dM \\
        & \ge \delta^{m-1} \sum_{i}r(\gamma(t_i))\\
        &\ge \delta^{m-1}\frac 14R.
    \end{split}
\end{equation*}
The result is proved.\\
\cqd

For complete surfaces in $\real^n$ that satisfy the Gauss-Bonnet relation, we obtain the following result.

\begin{mycoro} Let $\delta$ be as in theorem \ref{thmmean}. If $M$ is a complete noncompact surface in $\real^n$ satisfying
\begin{equation}2\pi\chi (M)-\int_MKdM<2\delta,\label{GB}\end{equation}
where $\chi(M)$ is the Euler characteristic of $M$, then
$$\int_M|H|dM=+\infty.$$\label{coro4.1}
\end{mycoro}
{\bf Proof}. From Theorem A of Shiohama \cite{Sh}, for any $q\in M$,
$$\lim_{r\rightarrow \infty}\frac{2V(B_r(q))}{ r^2}=2\pi\chi(M)-\int_MKdM.$$
Observe that there is a misprint in the denominator of this expression in Shiohama's paper. So,
$$\lim_{r\rightarrow \infty}\frac{V(B_r(q)}{\pi r^2}<\delta.$$
Thus, Theorem \ref{thmmean} implies the result.\\
\cqd
\begin{myrem}
The flat plane embedded in $\real^n$ shows that the condition (\ref{GB}) is necessary.
\end{myrem}

\parbox[t]{3.3in}{
Hil\'ario Alencar \\
Insitituto de Matem\'{a}tica\\
Universidade Federal de Alagoas\\
57072-900 Macei\'o-AL, Brazil\\
hilario@mat.ufal.br}
\parbox[t]{3.3in}{
Walcy Santos \\
Instituto de Matem\'{a}tica\\
Universidade Federal do Rio de Janeiro\\
Caixa Postal 68530\\
21941-909, Rio de Janeiro-RJ,
Brazil\\
walcy@im.ufrj.br}
\parbox[t]{3.3in}{
 Detang Zhou\\ Insitituto de Matem\' atica\\
Universidade Federal Fluminense\\   24020-140, Niter\'{o}i-RJ Brazil \\
zhou@impa.br}
 \end{document}